\documentclass{amsart}
\advance\headheight by 1.4778pt 
\usepackage{euscript,url}

\theoremstyle{plain}
\newtheorem*{corollary*}{Corollary}
\newtheorem{proposition}{Proposition}
\newtheorem{theorem}{Theorem}
\newtheorem{lemma}{Lemma}

\theoremstyle{definition}
\newtheorem*{definition*}{Definition}

\theoremstyle{remark}
\newtheorem*{remark}{Remark}
\newtheorem{example}{Example}

\let\epsilon\varepsilon
\def\sprod#1#2{\left\langle#1,#2\right\rangle}

\def\wslim{\mbox{w$^*$-}\mathop{\rm lim}\nolimits}
\def\colon{\;:\;}
\def\txt#1{\quad\mbox{#1}\quad}
\def\g{\mbox{\boldmath$\gamma$}}
\def\t{\mbox{\boldmath$\tau$}}
\def\a{\mbox{\boldmath$\alpha$}}
\def\os#1{{\EuScript #1}}

\def\N{\mathbb N{\mathstrut}}
\def\L{\mathbb L{\mathstrut}}
\def\F{\mathbb F{\mathstrut}}
\def\D{\mathbb D{\mathstrut}}

\begin{document}
\title[Uniformly convex operators and martingale type]
{Uniformly convex operators\\
  and martingale type}
\author{J\"org Wenzel}
\date{August 31, 1998}
\address{Mathematical Institute, Friedrich-Schiller-University Jena,
  07740 Jena, Germany}
\email{wenzel@minet.uni-jena.de}
\begin{abstract}
  The concept of uniform convexity of a Banach space was generalized
  to linear operators between Banach spaces and studied by
  Beauzamy~\cite{beauzamy76:_op}. Under this generalization, a Banach
  space $X$ is uniformly convex if and only if its identity map $I_X$
  is. Pisier showed that uniformly convex Banach spaces have
  martingale type $p$ for some $p>1$. We show that this fact is in
  general not true for linear operators. To remedy the situation, we
  introduce the new concept of martingale subtype and show, that it is
  equivalent, also in the operator case, to the existence of an
  equivalent uniformly convex norm on~$X$. In the case of identity
  maps it is also equivalent to having martingale type~$p$ for some
  $p>1$.

  Our main method is to use sequences of ideal norms defined on the
  class of all linear operators and to study the factorization of the
  finite summation operators. There is a certain analogy with the
  theory of Rademacher type.
\end{abstract}

\maketitle

\section{Introduction}
\label{sec:intro}

Banach spaces admitting an equivalent uniformly convex norm enjoy
several equivalent characterizations. Among others, they are the
superreflexive Banach spaces, i.~e. not only is such a Banach space
$X$ reflexive, but every Banach space whose finite dimensional
subspaces can be found (uniformly) in $X$ is reflexive.

A connection with martingales was studied by Pisier \cite{pis75}. For
$1<p\leq2$, a Banach space $X$ has martingale type $p$ if there exists
a constant $c\geq 0$ such that
\[ \Big(\int_0^1 \Big\|\sum_{k=1}^n d_k(t) \Big\|^p dt \Big)^{1/p}
\leq c \Big(\sum_{k=1}^n \int_0^1 \|d_k(t)\|^p dt \Big)^{1/p}
\]
for all $X$-valued martingale difference sequences $d_1,\dots,d_n$.

For our purpose the fundamental result of Pisier's paper
\cite[Thm.~3.2, p.~340]{pis75} combined with James's and Enflo's
investigations \cite[Thm.~4, p.~903]{james72:_super_banac},
\cite[p.~281]{enf72} can be summarized as follows.  (See further down
for detailed definitions.)
\begin{theorem}
  \label{theorem:2}
  For a Banach space $X$ the following properties are equivalent:
  \begin{enumerate}
    \itemsep=-.5\itemsep
  \item \label{item:3} $X$ has martingale type $p$ for some $p>1$,
  \item $X$ is superreflexive,
  \item $X$ admits an equivalent uniformly convex norm,
  \item $X$ does not factor the finite summation operators $\Sigma_n$
    uniformly.
  \end{enumerate}
\end{theorem}

If one is interested in linear operators between Banach spaces instead
of just Banach spaces, the above theorem is no longer true, as can be
seen by examples (see Example~\ref{ex:2} at the end).  The aim of this
paper is to prove an operator-theoretic replacement of
Theorem~\ref{theorem:2}; see Theorem~\ref{theorem:1} in
Section~\ref{sec:Conclusion}.

The main difficulties arising in this setting are the lack of a
suitable substitute of J-convexity in the operator case (see Beauzamy
\cite[p.~265]{bea85a} for a definition of J-convexity) and the fact
that the submultiplicativity of the martingale type ideal norms can no
longer be exploited.

The following general method has turned out to be useful to generalize
Banach space-theoretic results to results about operators.

Given a sequence of parameters $(\a_n)$ associating with every
operator $T:X\to Y$ a sequence of non-negative numbers $(\a_n(T))$,
let $\a_n(\infty):=\sup\a_n(T)$, where the supremum is taken over all
operators $T:X\to Y$ of norm $1$ and all Banach spaces $X$ and $Y$.
Then the sequence $(\a_n(\infty))$ describes the `worst' behavior that
can occur for an operator $T$. In the Banach space case, one is mostly
interested in the boundedness of the sequence $(\a_n(I_X))$ of the
identity map of a Banach space $X$, i.~e. $\a_n(I_X)=O(1)$.  In the
operator case, the behavior $\a_n(T)/\a_n(\infty)\to0$, i.~e.
$\a_n(T)=o(\a_n(\infty))$ is much more useful. Another beautiful
example for this heuristic in the context of Rademacher and Gauss type
is given by Hinrichs in \cite{hinrichs:_operat_radem_gauss}.

In particular, taking as $\a_n(T)$ the martingale type ideal norm
$\t_n(T)$ formed with $n$ martingale differences, the condition
$\t_n(T)=o(\sqrt n)$ will be the right replacement for
Condition~\ref{item:3} in Theorem~\ref{theorem:2}.

Let us quickly review the contents of this article. In
Section~\ref{sec:martingale} we introduce martingale and Haar (co)type
ideal norms, which are close relatives. The main result of this
section is that all the $o$-conditions described above for these four
types of ideal norms yield equivalent properties. In
Section~\ref{sec:equal} we establish the connection of the martingale
type ideal norms with the factorization of the finite summation
operators. To do so, we use a variant of the martingale type ideal
norms, namely the equal norm martingale type ideal norms. In
Section~\ref{sec:convexity} we repeat the definitions of uniform
convexity and uniform smoothness of linear operators, introduce the
super weakly compact operators and give various characterizations due
to Beauzamy \cite{beauzamy76:_op} connecting the two concepts.  Since
the main emphasis of this article is on the connection with
martingales, we are rather brief here and give mostly references for
the proofs. Finally in Section~\ref{sec:Conclusion} we formulate and
prove our main theorem and provide an example of an operator, for
which Theorem~\ref{theorem:2} is false.

To finish this introduction, let us point out some notational
conventions used throughout.  We write $B_X$ for the unit ball of a
Banach space $X$ and $I_X$ for its identity map. Furthermore, the
reader has already realized, that we use Landau's big-$O$ and
little-$o$ notation.

\section{Martingale type and cotype}
\label{sec:martingale}

First of all, we introduce the martingale type and cotype ideal norms.
They were first considered by Pisier in \cite[Rem.~3.3,
p.~346]{pis75}. We also consider Haar type and cotype ideal norms as
restrictions of martingale type and cotype ideal norms to special
classes of martingales. An operator is said to have the corresponding
subtype or subcotype, if these sequences of ideal norms behave just a
little better than the worst case. The main result of this section
will be that all four possible subtype and subcotype properties
coincide.

Let us start by considering an arbitrary martingale $(f_k)$ defined on
$[0,1)$ with values in a Banach space $X$ and adapted to a filtration
$\mathcal{F}_0 \subseteq \mathcal{F}_1 \subseteq \dots$ of finitely
generated $\sigma$-algebras. We denote by $E_k$ the operator of
conditional expectation with respect to the $\sigma$-algebra
$\mathcal{F}_k$.

In particular, taking for $\mathcal{F}_k$ the $\sigma$-algebra
generated by the dyadic intervals
$\Delta_k^{\!(j)}:=[\frac{j-1}{2^k},\frac j{2^k})$, we obtain the so
called \emph{Walsh-Paley} or \emph{dyadic} martingales. Since the
corresponding martingales $f_n$ are just linear combinations of the
Haar functions $\chi_k^{(j)}$ for $k=0,\dots,n$ and
$j=1,\dots,2^{k-1}$, we also use the term \emph{Haar polynomials} to
denote the functions $f_n$ in this case. Remember, that the Haar
functions $\chi_k^{(j)}$ are defined by
\[   \chi_k^{(j)}(t) :=
\begin{cases}
  +2^{(k-1)/2} & \txt{if $t\in\Delta_k^{\!(2j-1)}$,} \\
  -2^{(k-1)/2} & \txt{if $t\in\Delta_k^{\!(2j)}$,} \\
  0  & \txt{otherwise.}
\end{cases}
\]
We let $\chi_0^{(0)}\equiv 1$.

We will mainly deal with the sequence of martingale differences
$(d_k)$ instead of $(f_k)$, where $d_1=f_1$ and $d_k=f_k-f_{k-1}$ for
$k>1$. The Banach space of square integrable $X$-valued functions on
$[0,1)$ is denoted by $[L_2,X]$. For $f\in[L_2,X]$ we write
\[ \|f|L_2\| := \Big(\int_0^1 \|f(t)\|^2 \, dt\Big)^{1/2}.
\]
The fact that $\|E_kf|L_2\|\leq\|f|L_2\|$ for all $f\in[L_2,X]$ and
$k\in\N$ will be frequently used.

With each operator $T:X\to Y$ we associate the operator
$[L_2,T]:[L_2,X]\to [L_2,Y]$, defined by
\[ [L_2,T]f(t) := T\big(f(t)\big).
\]
For any two functions $f\in[L_2,X]$ and $g'\in[L_2,Y']$ we write
\[ \big\langle[L_2,T]f,g'\big\rangle :=
\int_0^1 \big\langle{Tf(t)},{g'(t)}\big\rangle\, dt.
\]

\begin{definition*}
  For $T:X\to Y$, the $n$-th \emph{martingale type ideal norm}
  $\t(T|\os M_n)$ is the smallest constant $c\geq 0$ such that
  \[ \Big\| \sum_{k=1}^n [L_2,T]d_k \Big| L_2 \Big\|
  \leq c \Big( \sum_{k=1}^n \|d_k|L_2\|^2 \Big)^{1/2}
  \]
  for all $X$-valued martingale difference sequences $d_1,\dots,d_n$
  adapted to any filtration on $[0,1)$.

  The $n$-th \emph{martingale cotype ideal norm}
  $\g(T|\os M_n)$ is the smallest constant $c\geq 0$ such that
  \[ \Big( \sum_{k=1}^n \big\|[L_2,T]d_k\big|L_2\big\|^2 \Big)^{1/2}
  \leq c \Big\| \sum_{k=1}^n d_k \Big| L_2 \Big\|
  \]
  for all $X$-valued martingale difference sequences $d_1,\dots,d_n$
  adapted to any filtration on $[0,1)$.
\end{definition*}
Note that, for norm one operators, both sequences $(\t(T|\os M_n))$
and $(\g(T|\os M_n))$ are bounded by $2\sqrt n$ and
therefore the following definition makes sense.
\begin{definition*}
  An operator $T:X\to Y$ has \emph{martingale subtype} or
  \emph{martingale subcotype} if $\t(T|\os M_n)/\sqrt n\to0$ or
  $\g(T|\os M_n)/\sqrt n\to0$, respectively.
\end{definition*}
For convenience we introduce the following notation for the
\emph{dyadic trees:}
\[ \D_m^n := \{ (k,j) \colon k=m,\dots,n;\ j=1,\dots,2^{k-1} \}
\txt{for $1\leq m \leq n$.}
\]
We let $\D_0^n := \D_1^n \cup \{ (0,0) \}$.
\begin{definition*}
  For $T:X\to Y$, the \emph{Haar type ideal norm} $\t(T|\os
  H(\D_m^n))$ associated with the index set $\D_m^n$ is the smallest
  constant $c\geq 0$ such that
  \[ \Big\| \sum_{\D_m^n} Tx_k^{(j)}\chi_k^{(j)}
  \Big| L_2 \Big\|
  \leq c
  \Big( \sum_{\D_m^n} \| x_k^{(j)} \|^2 \Big)^{1/2}
  \] for all $(x_k^{(j)})\subseteq X$.

  The \emph{Haar cotype ideal norm} $\g(T|\os H(\D_m^n))$
  associated with the index set $\D_m^n$ is the smallest
  constant $c\geq 0$ such that
  \[ \Big( \sum_{\D_m^n} \| Tx_k^{(j)} \|^2 \Big)^{1/2}
  \leq c
  \Big\| \sum_{\D_m^n} x_k^{(j)}\chi_k^{(j)}
  \Big| L_2 \Big\|
  \] for all $(x_k^{(j)})\subseteq X$.
\end{definition*}
As before, we will be interested in the suboptimal behavior of these
sequences.
\begin{definition*}
  An operator $T:X\to Y$ has \emph{Haar subtype} or \emph{Haar
    subcotype} if $\t(T|\os H(\D_0^n))/\sqrt n\to0$ or $\g(T|\os
  H(\D_0^n))/\sqrt n\to0$, respectively.
\end{definition*}

It is easily verified that
\begin{eqnarray*}
  &\t(T|\os H(\D_1^n))
  \leq \t(T|\os H(\D_0^n))
  \leq 2 \, \t(T|\os H(\D_1^n)),& \\
  &\g(T|\os H(\D_1^n))
  \leq \g(T|\os H(\D_0^n))
  \leq 3 \, \g(T|\os H(\D_1^n)).&
\end{eqnarray*}

Since for any $(x_k^{(j)})\subseteq X$ the sequence of functions $ d_k
:= \sum_{j=1}^{2^{k-1}} x_k^{(j)} \chi_k^{(j)} $ forms a sequence of
martingale differences, we obviously have
\begin{equation}\label{eq:1}
  \t(T|\os H(\D_1^n)) \leq \t(T|\os M_n) \txt{and}
  \g(T|\os H(\D_1^n)) \leq \g(T|\os M_n).
\end{equation}
No reverse estimate is known. In the limiting case where one of these
sequences of ideal norms behaves like $o(\sqrt n)$ we have, however,
equivalence.
\begin{proposition}\label{p:4}
  For any operator $T:X\to Y$ the following properties are equivalent:
  \begin{enumerate}
    \itemsep=-.5\itemsep
  \item $T$ has martingale subtype,
  \item $T$ has martingale subcotype,
  \item $T$ has Haar subtype,
  \item $T$ has Haar subcotype.
  \end{enumerate}
\end{proposition}

We postpone the proof of Proposition~\ref{p:4} in order to provide
some prerequisites. The main idea of the proof is contained in the
following lemma.
\begin{lemma}
  \label{l:8}
  Let $T:X\to Y$, then we have
  \[ \frac{\t(T|\os M_{2^n})}{2^{n/2}} \leq 3\, \frac{\g(T|\os
    H(\D_0^n))}{n^{1/2}}.
  \]
\end{lemma}
\begin{proof}
Given any sequence $(d_1,\dots,d_{2^n})$ of $X$-valued martingale
differences, define a function $F_n:[0,1)\times[0,1)\to X$ by
$F_n(s,t):=d_i(t)$ if $s\in\Delta_n^{\!(i)}$.

Obviously
\[ \int_0^1\!\!\int_0^1\|F_n(s,t)\|^2\, dt\,ds = \sum_{i=1}^{2^n}
\int_{\Delta_n^{\!(i)}} \!\!\int_0^1 \|F_n(s,t)\|^2\,dt\,ds =
\frac1{2^n} \sum_{i=1}^{2^n} \|d_i|L_2\|^2.
\]
Viewing $F_n(\,\cdot\,,t)$ as a Haar polynomial, we have
\[ F_n(s,t) = \sum_{\D_0^n} x_k^{(j)}(t)\,\chi_k^{(j)}(s)
\]
where the functions $x_k^{(j)}$ are defined by $x_k^{(j)}(t):=
\int\nolimits_0^1F_n(s,t)\,\chi_k^{(j)}(s)\,ds$. By the definition of
$\g(T|\os H(\D_0^n))$ we have
\[ \sum_{\D_0^n} \|Tx_k^{(j)}(t)\|^2 \leq \g(T|\os H(\D_0^n))^2\,
\Big\|\sum_{\D_0^n} x_k^{(j)}(t)\,\chi_k^{(j)}
\Big|L_2\Big\|^2.
\]
Integration with respect to $t\in[0,1)$ yields
\begin{eqnarray}
  \sum_{\D_0^n} \|[L_2,T]x_k^{(j)}|L_2\|^2
  & \leq &
  \g(T|\os H(\D_0^n))^2\, \|F_n|L_2\|^2 \nonumber\\
  \label{eq:8}
  & = &
  \g(T|\os H(\D_0^n))^2 \,\frac1{2^n}\sum_{i=1}^{2^n} \|d_i|L_2\|^2.
\end{eqnarray}

Letting $\N_k^{(j)} := \{i\colon\Delta_n^{\!(i)} \subseteq
\Delta_k^{\!(j)}\}$ and observing that $\N_{k-1}^{(j)} =
\N_k^{(2j-1)} \cup \N_k^{(2j)}$, we obtain for $x_k^{(j)}$
\begin{eqnarray*}
  x_k^{(j)}(t)
  & = &
  \int_0^1 F_n(s,t)\,\chi_k^{(j)}(s)\,ds
  =
  2^{(k-1)/2}\, \Big(\!\!\! \int_{\Delta_k^{\!(2j-1)}} \!\!\!
  F_n(s,t)\,ds - \!\!\! \int_{\Delta_k^{\!(2j)}}\!\! F_n(s,t)\,
  ds\Big) \\
  & = &
  2^{(k-1)/2-n}\, \Big( \!\! \sum_{\N_k^{(2j-1)}} \!\! d_i(t) -
  \!\! \sum_{\N_k^{(2j)}} \!d_i(t) \Big).
\end{eqnarray*}
This implies that
\begin{equation}
  \label{eq:9}
  \|[L_2,T]x_k^{(j)}|L_2\| = 2^{(k-1)/2-n} \, \Big\|
  \!\!\sum_{\N_k^{(2j-1)}} \!\! [L_2,T]d_i - \!\sum_{\N_k^{(2j)}}\!
  [L_2,T]d_i \Big| L_2 \Big\|.
\end{equation}
Since the conditional expectation operator has norm one in $[L_2,X]$
it now follows that
\[ \Big\|\!\!\sum_{\N_k^{(2j-1)}} \!\! [L_2,T]d_i\Big|L_2\Big\|
  \leq \Big\|\!\!\sum_{\N_k^{(2j-1)}}\! \! [L_2,T]d_i - \!
  \sum_{\N_k^{(2j)}} \! [L_2,T]d_i\Big|L_2\Big\|
\]
and therefore by the triangle inequality
\begin{eqnarray}
  \Big\|\!\sum_{\N_{k-1}^{(j)}} \![L_2,T]d_i \Big|
  L_2\Big\|
  & \leq &
  \Big\|\!\!\sum_{\N_k^{(2j-1)}} \!\! [L_2,T]d_i - \!
  \sum_{\N_k^{(2j)}} \! [L_2,T]d_i\Big|L_2\Big\|
  +
  2\, \Big\|\!\!\sum_{\N_k^{(2j-1)}} \!\!
  [L_2,T]d_i\Big|L_2\Big\| \nonumber \\
  & \leq &
  3\, \Big\|\!\!\sum_{\N_k^{(2j-1)}}\! \! [L_2,T]d_i - \!
  \sum_{\N_k^{(2j)}} \! [L_2,T]d_i\Big|L_2\Big\|.
  \label{eq:10}
\end{eqnarray}
Using~(\ref{eq:8}), (\ref{eq:9}) and~(\ref{eq:10}) the proof completes
as follows:
\begin{eqnarray*}
\lefteqn{
    \Big\| \sum_{i=1}^{2^n} [L_2,T]d_i \Big| L_2 \Big\|
    =
    \Bigg( \frac1n \sum_{k=1}^n
    \Big\| \sum_{j=1}^{2^{k-1}} \sum_{\N_{k-1}^{(j)}}
    [L_2,T]d_i \Big| L_2 \Big\|^2 \Bigg)^{1/2}} \\
    & \leq &
    \Bigg( \frac1n \sum_{k=1}^n
    \Big( \sum_{j=1}^{2^{k-1}} \Big\| \sum_{\N_{k-1}^{(j)}}
    [L_2,T]d_i \Big| L_2 \Big\| \Big)^2 \Bigg)^{1/2} \\
    & \leq &
    3\, \Bigg( \frac1n \sum_{k=1}^n 2^{k-1}
    \sum_{j=1}^{2^{k-1}} \Big\| \sum_{\N_k^{(2j-1)}}
    [L_2,T]d_i -\sum_{\N_k^{(2j)}} [L_2,T]d_i \Big| L_2
    \Big\|^2 \Bigg)^{1/2} \\
    & \leq &
    3 \,\Big( \frac{2^{2n}}n \sum_{k=1}^n \sum_{j=1}^{2^{k-1}} \|
    [L_2,T]x_k^{(j)} | L_2 \|^2 \Big)^{1/2}
    \, \leq \,
    3 \,\g(T|\os H(\D_0^n))\, \Big( \frac{2^n}n \sum_{i=1}^{2^n}
    \|d_i|L_2\|^2 \Big)^{1/2}.
\end{eqnarray*}
\end{proof}

We next observe that the martingale and Haar type and cotype ideal
norms are dual to each other.
\begin{proposition}
  \label{l:9}
  For $0\leq m<n$ we have
  \begin{eqnarray*}
\arraycolsep=0pt
    \g(T\phantom{'}|\os H(\D_m^n)) \leq & \t(T'|\os H(\D_m^n))
    & \leq 2\,\g(T\phantom{'}|\os H(\D_m^n)), \\
    \g(T'|\os H(\D_m^n)) \leq & \t(T\phantom{'}|\os H(\D_m^n))
    & \leq 2\,\g(T'|\os H(\D_m^n)).
  \end{eqnarray*}
  If $m=0$, we can omit the factors $2$ and have equality. On the
  other hand
  \begin{eqnarray*}
\arraycolsep=0pt
    \g(T\phantom{'}|\os M_n) \leq & 2\,\t(T'|\os M_n) & \leq
    4\,\g(T\phantom{'}|\os M_n), \\
    \g(T'|\os M_n) \leq & 2\,\t(T\phantom{'}|\os M_n) & \leq
    4\,\g(T'|\os M_n).
  \end{eqnarray*}
\end{proposition}
\begin{proof}
  The proof can be obtained using standard duality techniques and is
  left to the reader.
\end{proof}

We can now prove Proposition~\ref{p:4}.
\begin{proof}
It follows from Lemma~\ref{l:8} and~(\ref{eq:1}) that
\[ \frac{\t(T|\os M_{2^n})}{2^{n/2}} \leq 3\, \frac{\g(T|\os
  H(\D_0^n))}{n^{1/2}} \leq 6\, \frac{\g(T|\os M_n))}{n^{1/2}}.
\]
Since the same is true for $T'$, it follows from Proposition~\ref{l:9}
that
\[ \frac{\g(T|\os M_{2^n})}{2^{n/2}} \leq 6\, \frac{\t(T|\os
  H(\D_0^n))}{n^{1/2}} \leq 12\, \frac{\t(T|\os M_n))}{n^{1/2}}.
\]
Hence if one of these quotients tends to zero, all the others tend to
zero too, which proves the proposition by virtue of the monotonicity
of the involved ideal norms.
\end{proof}

\section{Equal norm martingale type}
\label{sec:equal}

\begin{definition*}
  For $T:X\to Y$, the $n$-th \emph{equal norm martingale type ideal
    norm} $\t^\circ(T|\os M_n)$ is the smallest constant $c\geq 0$
  such that
  \[ \Big\| \sum_{k=1}^n [L_2,T]d_k \Big| L_2 \Big\|
  \leq c \Big( \sum_{k=1}^n \|d_k|L_2\|^2 \Big)^{1/2}
  \]
  for all $X$-valued martingale difference sequences $d_1,\dots,d_n$
  adapted to any filtration on $[0,1)$ under the additional assumption
  that $\|d_1|L_2\|=\dots=\|d_n|L_2\|$.
\end{definition*}
The quantities $\t^\circ(T|\os M_n)$ can also be defined in a
different way.
\begin{lemma}
  \label{lem:4}
  For $T:X\to Y$, the ideal norm $\t^\circ(T|\os M_n)$ is the smallest
  constant $c\geq 0$ such that
  \begin{equation}
    \label{eq:13}
    \Big\| \sum_{k=1}^n [L_2,T]d_k \Big| L_2 \Big\|
    \leq c \,n^{1/2} \sup_{k=1,\dots,n} \|d_k|L_2\|
  \end{equation}
  for all $X$-valued martingale difference sequences $d_1,\dots,d_n$
  adapted to any filtration on $[0,1)$.
\end{lemma}
\begin{proof}
  For the time being denote by $\t^{\circ\circ}(T|\os M_n)$ the
  smallest constant such that~(\ref{eq:13}) holds. It is obvious that
  for $\|d_1|L_2\|=\dots=\|d_n|L_2\|$
\begin{eqnarray*}
  \Big\|\sum_{k=1}^n [L_2,T]d_k \Big|L_2\Big\|
  & \leq &
  \t^{\circ\circ}(T|\os M_n) \, n^{1/2} \, \sup_{k=1,\dots,n}
  \|d_k|L_2\|^2 \\
  & = &
  \t^{\circ\circ}(T|\os M_n) \, \Big(\sum_{k=1}^n
  \|d_k|L_2\|^2\Big)^{1/2}.
\end{eqnarray*}
Therefore $\t^\circ(T|\os M_n) \leq \t^{\circ\circ}(T|\os M_n)$.

On the other hand, let $d_1,\dots,d_n$ be an arbitrary sequence of
martingale differences. It follows that for $\tilde d_k :=
d_k/\|d_k|L_2\|$
\[ \Big\|\sum_{k=1}^n [L_2,T]\tilde d_k\Big|L_2\Big\| \leq
\t^{\circ}(T|\os M_n) \, n^{1/2}.
\]
But the same is true for $\zeta_k\tilde d_k$ instead of $\tilde d_k$,
where $|\zeta_k|=1$. An extreme point argument then yields that
\[ \Big\|\sum_{k=1}^n [L_2,T]\alpha_k \tilde d_k\Big|L_2\Big\| \leq
\t^{\circ}(T|\os M_n) \, n^{1/2}
\]
whenever $|\alpha_k|\leq 1$. In particular, we may take
\[ \alpha_k:=\frac{\|d_k|L_2\|}{\sup\limits_{h=1,\dots,n}\|d_h|L_2\|},
\]
which shows that $\t^{\circ\circ}(T|\os M_n) \leq \t^{\circ}(T|\os
M_n)$.
\end{proof}

Obviously we have
\[ \t^\circ(T|\os M_n) \leq \t(T|\os M_n).
\]
The main purpose of this section is to prove a reverse estimate. We
follow an approach similar to Bourgain\slash Kalton\slash Tzafriri in
\cite[Thm.~3.1., p.~160]{bourgain89:_geomet_l} where they show that
equal norm Rademacher type $2$ is equivalent to ordinary Rademacher
type $2$. The main idea is contained in the following construction of
`glueing' together $m$ copies of a given martingale of length $n$,
which yields a martingale of length $mn$ with smaller differences. An
appropriate blocking of this longer martingale will then give a
martingale of length of order $n$ and nearly equal $L_2$-norms.

Let $\phi_j^m:[\frac{j-1}m,\frac jm)\to [0,1)$ be defined by
$\phi_j^m(t):=mt-j+1$. Given a function $f:[0,1)\to X$, we denote by
$\Phi_j^mf$ the function
\[ \Phi_j^mf(t) :=
\begin{cases}
  f(\phi_j^m(t)) & \mbox{if $t\in[\frac{j-1}m,\frac jm)$,} \\
  0 & \mbox{otherwise.}
\end{cases}
\]
Given a sequence of martingale differences $(\sigma_k)$,
the sequence of martingale differences
\[ \Phi_1^md_1,\dots,\Phi_m^md_1, \Phi_1^md_2,\dots,\Phi_m^md_2,
\dots, \Phi_1^md_n,\dots,\Phi_m^md_n
\]
is adapted to the filtration
\[ \Phi_1^m\mathcal{F}_1,\dots,\Phi_m^m\mathcal{F}_1,
\Phi_1^m\mathcal{F}_2,\dots,\Phi_m^m\mathcal{F}_2, \dots,
\Phi_1^m\mathcal{F}_n,\dots,\Phi_m^m\mathcal{F}_n,
\]
where $\Phi_j^m\mathcal{F}_k$ is the $\sigma$-algebra generated by all
sets $A\subseteq[0,1)$ such that $\phi_j^m(A)\in\mathcal{F}_k$ and by
all its predecessor $\sigma$-algebras.

Observe that
\begin{equation}
  \label{eq:14}
  \Big\| \sum_{k=1}^n \sum_{j=1}^m \Phi_j^md_k \Big|L_2\Big\| =
  \Big\|\sum_{k=1}^n d_k\Big| L_2\Big\|
  \txt{and}
  \|\Phi_j^md_k|L_2\|^2 = \frac1m \, \|d_k|L_2\|^2.
\end{equation}
Moreover, all differences $\Phi_j^md_k$ in any block of length at most
$m$ have disjoint support.

\begin{lemma}
  \label{lem:6}
  The sequence $(\t^\circ(T|\os M_n))$ is non-decreasing.
\end{lemma}
\begin{proof}
  Let $d_1,\dots,d_n$ be $X$-valued martingale differences such that
  $\|d_k|L_2\|=1$. For $m:=n+1$ the construction above yields a
  martingale difference sequence of length $n(n+1)$. Define a new
  sequence of martingale differences $\tilde d_1,\dots,\tilde d_{n+1}$
  by blocking $n$ consecutive terms:
\[ \tilde d_h := \Phi_{n-h+3}^md_{h-1} + \dots + \Phi_{n+1}^md_{h-1} +
\Phi_1^md_h + \dots \Phi_{n-h+1}^md_h.
\]
Since $\|\tilde d_h|L_2\| = \sqrt{\frac n{n+1}}$ it follows that
\[ \Big\|\sum_{k=1}^n [L_2,T]d_k \Big| L_2\Big\| =
\Big\|\sum_{h=1}^{n+1} [L_2,T]\tilde d_h \Big|L_2\Big\| \leq
\t^\circ(T|\os M_{n+1})\, \Big(\sum_{k=1}^n \|d_k|L_2\|^2\Big)^{1/2},
\]
which proves that $\t^\circ(T|\os M_n) \leq \t^\circ(T|\os M_{n+1})$.
\end{proof}

\begin{lemma}
  \label{lem:3}
  $\quad \t(T|\os M_n) \leq 16\,\t^\circ(T|\os M_n)$.
\end{lemma}
\begin{proof}
  Let $d_1,\dots,d_n$ be $X$-valued martingale differences. By scaling
  we may assume that $\sum_{k=1}^n \|d_k|L_2\|^2 = 1$.

Let $l$ be such that $4^l\leq 16\,n < 4^{l+1}$ and $m:=4^l$. For
$h=1,2\dots$ define
\[ \F_h := \Big\{ k \colon \frac1{2^h} < \|d_k|L_2\| \leq
\frac2{2^h} \Big\}
\txt{and}
\F  := \bigcup_{h=1}^l \F_h.
\]
First of all, we estimate the sum of all the differences with small
norm:
\begin{eqnarray}
  \Big\|\sum_{k\not\in\F } [L_2,T]d_k \Big|L_2\Big\|
  & \leq &
  \t(T|\os M_n)\, \Big(\sum_{k\not\in\F } \|d_k|L_2\|^2
  \Big)^{1/2} \nonumber \\
  & \leq &
  \t(T|\os M_n)\, (n4^{-l})^{1/2}
  \leq
  \frac 12 \t(T|\os M_n).\label{eq:7}
\end{eqnarray}
Here we used, that the sequence $(\t(T|\os M_n))$ is obviously
non-decreasing and that $16n<4^{l+1}$.

For the martingale difference sequence $(d_k)$ with $k\in\F $ we
apply the glueing technique described above. Then it follows
from~(\ref{eq:14}) that for $k\in\F_h$
\[ \frac1{m4^h} < \|\Phi_j^md_k|L_2\|^2 = \frac1m\|d_k|L_2\|^2 \leq
\frac4{m4^h}.
\]
Therefore by disjointness, for any subset $\L\subseteq\{1,\dots,m\}$
of cardinality $|\L |=4^h$
\[ \frac1m < \Big\| \sum_{i\in\L } \Phi_i^md_k \Big|L_2\Big\|^2
\leq \frac4m.
\]
Writing $\L_j^{(h)} := \{4^h\,(j-1)+1, \dots 4^hj\}$ and
$\tilde d_j^k := \sum_{i\in\L_j^{(h)}} \Phi_i^md_k$ we obtain a
martingale of length
\[ N:= \sum_{h=1}^l \sum_{k\in\F_h} 4^{l-h} \leq 4^l
\sum_{h=1}^l \sum_{k\in\F_h} \|d_k|L_2\|^2 \leq 4^l \leq 16n.
\]
It follows from Lemma~\ref{lem:6} that $\t^\circ(T|\os
M_N)\leq\t^\circ(T|\os M_{16n})$ and therefore, we obtain from
Lemma~\ref{lem:4} that
\begin{eqnarray*}
  \Big\| \sum_{h=1}^l \sum_{k\in\F_h} \sum_{j=1}^{4^{l-h}}
  [L_2,T]\tilde d_j^k \Big| L_2\Big\|
  & \leq &
  \t^\circ(T|\os M_N)\, \sup \|\tilde d_j^k|L_2\| \sqrt N \\
  & \leq &
  \t^\circ(T|\os M_{16n})\,\frac2{\sqrt m} 2^l
  =
  2\,\t^\circ(T|\os M_{16n}).
\end{eqnarray*}
This shows that
\begin{equation}
  \label{eq:6}
  \Big\| \sum_{k\in\F } [L_2,T]d_k \Big| L_2\Big\| \leq
  2\,\t^\circ(T|\os M_{16n}).
\end{equation}
Putting together~(\ref{eq:7}) and~(\ref{eq:6}) we obtain
\[ \t(T|\os M_n) \leq \frac12 \t(T|\os M_n) + 2\,\t^\circ(T|\os
M_{16n}),
\]
which implies $\t(T|\os M_n) \leq 4\,\t^\circ(T|\os M_{16n})$. Finally
the assertion follows from the trivial fact that $\t^\circ(T|\os
M_{16n}) \leq 4\,\t^\circ(T|\os M_n)$.
\end{proof}

The \emph{summation operator} $\Sigma:l_1\to l_\infty$ is defined by
\[ \Sigma(\xi_k) := \Big(\sum_{h=1}^k \xi_h\Big), \]
while the \emph{finite summation operators} $\Sigma_n:l_1^n\to
l_\infty^n$ act between the finite dimensional spaces and are defined
in the same way.

The significance of the ideal norms $\t^\circ(T|\os M_n)$ is due to
the following fact, which establishes the connection with the
factorization of the finite summation operators.
\begin{proposition}
  \label{prop:3}
  There exists a factorization $\Sigma_n=B_n[L_2,T]A_n$ of the finite
  summation operator $\Sigma_n$, such that $\|B_n\|\,\|A_n\|\leq
  6\sqrt n/\t^\circ(T|\os M_{2n})$.
\end{proposition}
\begin{proof}
  There is nothing to prove for $T=0$. For $T\not=0$, by definition,
  for all $0<\delta<1$, there exists a sequence of martingale
  differences $d_1,\dots,d_{2n}$ such that $\|d_k|L_2\|=1$ and
  \[ \Big\|\sum_{k=1}^{2n} [L_2,T]d_k\Big| L_2\Big\| >
  \delta\,\t^\circ(T|\os M_{2n})\sqrt{2n}.
  \]
Choose $g'\in[L_2,Y']$ such that $\|g'|L_2\|=1$ and
\[ \Big\langle \sum_{k=1}^{2n} [L_2,T]d_k , g'\Big\rangle
> \delta\,\t^\circ(T|\os M_{2n})\sqrt{2n}.
\]
Let
\[ \F  := \Big\{ k\colon \big\langle{[L_2,T]d_k},{g'}\big\rangle >
\t^\circ(T|\os M_{2n}) \frac\delta{4\sqrt{2n}} \Big\}.
\]
It follows from Lemma~\ref{lem:6} that
\[ \Big\langle \sum_{k\in\F } [L_2,T]d_k, g'
\Big\rangle \leq \Big\|\sum_{k\in\F } [L_2,T]d_k\Big|L_2\Big\|
\,\|g'|L_2\| \leq \t^\circ(T|\os M_{2n}) |\F |^{1/2}
\]
and therefore
\begin{eqnarray*}
  \delta\,\t^\circ(T|\os M_{2n}) \sqrt{2n}
  & < &
  \Big\langle \sum_{k\in\F } [L_2,T]d_k , g'\Big\rangle
  + \sum_{k\not\in\F } \big\langle{[L_2,T]d_k},{g'}\big\rangle
  \\
  & \leq &
  \t^\circ(T|\os M_{2n}) |\F |^{1/2} + 2n\, \t^\circ(T|\os
  M_{2n}) \frac\delta{4\sqrt{2n}}.
\end{eqnarray*}
This shows that $m:=|\F | \geq \delta^2\,9n/8$. Choosing
$\delta$ appropriately, we may arrange that $m\geq n$. In particular,
we find elements $i_1<\dots<i_n$ in $\F $.

We can now define $A_n:l_1^n\to [L_2,X]$ by
\[ A_ne_k := \frac{d_{i_k}}{\big\langle{[L_2,T]d_{i_k}},{g'}\big\rangle}
\]
and $B_n:[L_2,Y] \to l_\infty^n$ by
\[ B_nf := \big( \sprod {f}{E_{i_k} g'}\big)_{k=1}^n.
\]
\end{proof}

\section{Uniform convexity and smoothness and super weakly compact
  operators}
\label{sec:convexity}

We now show, how the concepts above connect to the theory of weakly
compact operators.

First of all, we repeat the classical definitions of uniform convexity
and uniform smoothness of Banach spaces in the more general case of
linear operators between Banach spaces. See e.~g. Beauzamy
\cite[Def.~7, p.~121]{beauzamy76:_op}.
\begin{definition*}
  An operator $T:X\to Y$ is \emph{uniformly convex} if for all
  $\epsilon>0$ there exists $\delta>0$ such that for $\|x_{\pm}\|=1$
  with $\frac{\|x_++x_-\|}2\geq 1-\delta$ it follows that
  $\frac{\|Tx_+-Tx_-\|}2\leq \epsilon$.

  The operator $T$ is \emph{uniformly smooth} if for all $\epsilon>0$
  there exists $\delta>0$ such that for $\|y\|=1$ and
  $\|x\|\leq\delta$ it follows that
  $\frac{\|y+Tx\|+\|y-Tx\|}2\leq1+\epsilon\,\|x\|$.
\end{definition*}

Using an equivalent norm on the target space $Y$ does not spoil
uniform convexity. This is, however, not the case for the source space
$X$. Similarly, uniform smoothness depends on the special choice of
the norm on the target space.  We are therefore rather interested in
uniformly convex renormable and uniformly smooth renormable operators.
\begin{definition*}
  An operator $T:X\to Y$ is \emph{uniformly convex renormable} if
  there exists an equivalent norm $|||\,\cdot\,|||$ on $X$ such that
  $T:[X,|||\,\cdot\,|||]\to Y$ is uniformly convex.

  An operator $T:X\to Y$ is \emph{uniformly smooth renormable} if
  there exists an equivalent norm $|||\,\cdot\,|||$ on $Y$ such that
  $T:X\to [Y,|||\,\cdot\,|||]$ is uniformly smooth.
\end{definition*}
The above properties are equivalent to $T$ factoring through a
uniformly convex or uniformly smooth Banach space, respectively; see
\cite[Prop.~7.10.11]{pietsch:_orthon_banac}.

As in the case of Banach spaces, we can easily show that the concepts
of uniform convexity and smoothness are dual to each other. See
Lindenstrauss \cite[Thm.~1]{lindenstrauss63:_banac} for a proof in the
Banach space case, which can easily be carried over to linear
operators.
\begin{proposition}
  \label{p-convexity-dual}
  An operator $T:X\to Y$ is uniformly convex if and only if its dual
  $T':Y'\to X'$ is uniformly smooth.

  The operator $T$ is uniformly smooth if and only if $T'$ is
  uniformly convex.
\end{proposition}

Next we introduce super weakly compact operators, whose Banach space
counterparts are the superreflexive Banach spaces, and show their
connection with the factorization of the summation operators. Again,
this basically follows from the, by now, classical proof that
superreflexive Banach spaces do not factor the summation operators
uniformly, cf. James \cite{james72:_super_banac}, Beauzamy
\cite[Prop.~7, p.~236]{bea85a}, or Heinrich \cite[Thm.~5.1,
p.~29]{heinrich80:_finit}.

\begin{definition*}
  An operator $T:X\to Y$ is \emph{weakly compact}, if the image of the
  closed unit ball of $X$ under $T$ is relatively weakly compact in
  $Y$.
\end{definition*}

For the theory of ultraproducts of Banach spaces and linear operators,
we refer to Heinrich's papers
\cite{heinrich80:_finit,heinrich80:_ultrap_banac}. We mainly use the
notation of Pietsch\slash Wenzel \cite{pietsch:_orthon_banac}.

\begin{definition*}
  An operator $T:X\to Y$ is \emph{super weakly compact}, if all its
  ultrapowers $T^\mathcal{U}:X^\mathcal{U}\to Y^\mathcal{U}$ are
  weakly compact.
\end{definition*}

The summation operator is the typical non-weakly compact operator.
\begin{example}
  \label{ex:1}
  The summation operator $\Sigma$ is not weakly compact.
\end{example}

\begin{definition*}
  An operator $T:X\to Y$ is said to \emph{factor the summation
    operator} $\Sigma$, if there exist operators $A:l_1\to X$ and
  $B:Y\to l_\infty$ such that $\Sigma=BT\!A$.

  The operator $T$ is said to \emph{factor the finite summation
    operators $\Sigma_n$ uniformly}, if there exists a constant $c>0$,
  such that for all $n\in\N $ we can find factorizations
  $\Sigma_n=B_nT\!A_n$ such that $\|B_n\|\,\|A_n\|\leq c$.
\end{definition*}

The next proposition connects the above concepts.

\begin{proposition}
  \label{p-swc-summation-operator}
  Let $T:X\to Y$, then the following properties are equivalent:
  \begin{enumerate}
    \itemsep=-.5\itemsep
  \item\label{item:1} $T$ is super weakly compact,
  \item \label{item:4} $T^\mathcal{U}$ does not factor the summation
    operator $\Sigma$ for any ultrafilter $\mathcal{U}$,
  \item\label{item:6} $T$ does not factor the finite summation
    operators $\Sigma_n$ uniformly.
  \end{enumerate}
\end{proposition}
\begin{proof}
The equivalence of~\ref{item:1} and~\ref{item:4} is due to
Lindenstrauss/Pe\l czy\'nski \cite{lindenstrauss68:_absol_l}.

Assume that $T:X\to Y$ factors the finite summation operators
uniformly and let $\mathcal{U}$ be any non-trivial ultrafilter on
$\N$.

Let $J$ be the canonical map from $l_1$ into $l_1^\mathcal{U}$ induced
by the map $x\mapsto(x,x,\ldots)$. Note moreover, that the map
$(x_n)\mapsto \wslim_\mathcal{U} x_n$ induces a well defined operator
$Q$ from $l_\infty^\mathcal{U}$ onto $l_\infty$. The last fact is due
to the weak-$*$-compactness of the closed unit ball of $l_\infty$.
Obviously $\|J\|=\|Q\|=1$.

By the assumption, there exists a constant $c$, such that we find for
all $n\in\N $ operators $A_n:l_1^n\to X$ and $B_n:Y\to l_\infty^n$
satisfying $\|A_n\|\leq1$, $\|B_n\|\leq c$, and $\Sigma_n=B_nT\!A_n$.

Denote furthermore by $J_n$ the canonical embedding of $l_\infty^n$
into $l_\infty$ and by~$Q_n$ the projection from $l_1$ onto $l_1^n$.

We obtain a factorization of the summation operator $\Sigma$ via
$T^\mathcal{U}$ as
\[ \Sigma=Q(J_n\Sigma_nQ_n)^\mathcal{U}J
=Q(J_nB_n)^\mathcal{U}T^\mathcal{U}(A_nQ_n)^\mathcal{U}J.
\]
This proves that~\ref{item:4} $\Rightarrow$~\ref{item:6}.

If on the other hand~\ref{item:4} does not hold, then $T^\mathcal{U}$
factors $\Sigma$ for some ultrafilter $\mathcal{U}$. In the sense of
Heinrich \cite[Def.~1.1, p.~7]{heinrich80:_finit}, $T^\mathcal{U}$ is
finitely representable in $T$, i.~e.~for all $\epsilon>0$ and each
finite dimensional subspace $\widehat X_0\subseteq X^\mathcal{U}$ and
finite codimensional subspace $\widehat Y_0\subseteq Y^\mathcal{U}$
there are a finite dimensional subspace $X_0\subseteq X$ and a finite
codimensional subspace $Y_0\subseteq Y$ and maps $R:\widehat X_0\to
X_0$ and $S:Y/Y_0\to Y^\mathcal{U}/\widehat Y_0$ such that
\[ \big\|\widehat QT^\mathcal{U}\widehat J - SQTJR \big\| \leq \epsilon,
\]
where $J$ and $\widehat J$ are the canonical embedding maps of $X_0$
and $\widehat X_0$ and $Q$ and $\widehat Q$ are the canonical quotient
maps of $Y_0$ and $\widehat Y_0$.

It follows, that for all $n\in\N$ there is a subspace $X_0\subseteq X$
and a subspace $Y_0\subseteq Y$ and operators $R:l_1^n\to X_0$ and
$S:Y/Y_0\to l_\infty^n$, such that
\[ \big\|\Sigma_n - SQTJR\big\| \leq \epsilon.
\]
Moreover, since $\Sigma_n$ is injective, by the remark following
Definition 1.1 in \cite[p.~8]{heinrich80:_finit}, we can even arrange
that $\epsilon=0$.  This implies that $T$ factors $\Sigma_n$
uniformly.
\end{proof}

It was shown in Beauzamy \cite[Thm.~I.1, p.~111]{beauzamy76:_op} that
uniform convex renormability and super weak compactness are in fact
equivalent properties. See also Heinrich's paper
\cite[Thm.~5,p.~29]{heinrich80:_finit} and the detailed presentation
in Pietsch\slash Wenzel \cite[Sect.~7.6]{pietsch:_orthon_banac}.

One obtains the following equivalences.

\begin{proposition}
  \label{prop:5}
  Let $T:X\to Y$, then the following properties are equivalent:
  \begin{enumerate}
    \itemsep=-.5\itemsep
  \item $T$ is super weakly compact,\label{item:8}
  \item $T$ is uniformly convex renormable,\label{item:7}
  \item $T$ does not factor the finite summation operators $\Sigma_n$
    uniformly.\label{item:9}
  \end{enumerate}
\end{proposition}

\section{Main theorem}
\label{sec:Conclusion}

Before we can formulate and prove the main theorem some more
preparations are required.

The following classical result provides the missing link between the
operators $T$ and $[L_2,T]$. It was first proved by Day \cite[Thm.~2,
p.~504]{day41:_some} for Banach spaces, but Day's proof can
straightforwardly be extended to the operator case.

\begin{proposition}
  \label{p:3}
  The operator $T$ is uniformly convex if and only if $[L_2,T]$ is
  uniformly convex.
\end{proposition}

Finally, we will need the Haar cotype ideal norms of the finite
summation operators.

\begin{example}
  \label{ex:3}
  $ \quad\frac12\, (n+1)^{1/2} \leq \g_n(\Sigma_{2^n}|\os H(\D_0^n))
  \leq (n+1)^{1/2}.
  $
\end{example}
\begin{proof}
Obviously, for any operator $T: X\to Y$
\begin{eqnarray*}
  \Big\|\sum_{\D_0^n} Tx_k^{(j)}\chi_k^{(j)}\Big|L_2\Big\|
  & \leq &
  (n+1)^{1/2} \Big( \sum_{k=0}^n \Big\|\sum_{j=1}^{2^{k-1}}
  Tx_k^{(j)}\chi_k^{(j)} \Big|L_2\Big\|^2\Big)^{1/2} \\
  & = &
  (n+1)^{1/2} \Big( \sum_{k=0}^n \sum_{j=1}^{2^{k-1}}
  \|Tx_k^{(j)}\|^2 \|\chi_k^{(j)}|L_2\|^2\Big)^{1/2} \\
  & \leq &
  (n+1)^{1/2}\|T\|\, \Big(\sum_{\D_0^n}
  \|x_k^{(j)}\|^2\Big)^{1/2},
\end{eqnarray*}
where we have used that for fixed $k$ the Haar functions
$\chi_k^{(1)},\dots,\chi_k^{(2^{k-1})}$ have disjoint support.  This
proves the upper estimate by virtue of Proposition~\ref{l:9}.

To see the lower estimate, let $f:[0,1)\to l_1^{2^n}$ be defined by
\[ f(t) := e_i \quad\mbox{if $t\in\Delta_n^{\!(i)}$,}
\]
where $e_i$ is the $i$-th unit vector in $l_1^{2^n}$.  Write $f$ as a
Haar polynomial
\[ f=\sum_{\D_0^n} x_k^{(j)} \chi_k^{(j)},
\qquad\mbox{where}\quad x_k^{(j)}=\int_0^1 f(t)\chi_k^{(j)}(t)\, dt.
\]
Obviously
\[ \|f|L_2\|=\Big\|\sum_{\D_0^n} x_k^{(j)}\chi_k^{(j)}\Big|
L_2\Big\| = \Big( \frac1{2^n}\sum_{i=1}^{2^n} \|e_i\|^2\Big)^{1/2} =
1.
\]
Writing $\N_k^{(j)} := \{i \colon
\Delta_n^{\!(i)}\subseteq\Delta_k^{\!(j)} \}$, it follows that
\[ x_k^{(j)} = \frac1{2^n} \,2^{(k-1)/2}\,
(0,\dots,0\underbrace{+1,\dots,+1}_{\N_k^{(2j-1)}},
\underbrace{-1,\dots,-1}_{\N_k^{(2j)}},0,\dots,0)
\]
and
\[ \Sigma_{2^n}x_k^{(j)} = \frac1{2^n}\, 2^{(k-1)/2}\,
(0,\dots,0\underbrace{1,2,\dots,2^{n-k}}_{\N_k^{(2j-1)}},
\underbrace{2^{n-k}-1,\dots,1,0}_{\N_k^{(2j)}},0,\dots,0),
\]
which in turn yields $\|\Sigma_{2^n}x_k^{(j)}\| = 2^{-(k+1)/2}$. Hence
\[ \Big(\sum_{\D_0^n}\big\| \Sigma_{2^n}x_k^{(j)} \big\|^2
\Big)^{1/2} = \frac12 \, (n+1)^{1/2},
\]
which proves the assertion.
\end{proof}

We are now ready to state our most important result.
\begin{theorem}
  \label{theorem:1}
  For an operator $T:X\to Y$ the following properties are equivalent:
  \begin{enumerate}
    \itemsep=-.5\itemsep
  \item \label{item:15} $T$ has martingale subtype,
  \item \label{item:13} $T$ has martingale subcotype,
  \item \label{item:17} $T$ has Haar subtype,
  \item \label{item:18} $T$ has Haar subcotype,
  \item \label{item:11} $T$ is super weakly compact,
  \item \label{item:12} $T$ does not factor the finite summation
    operators $\Sigma_n$ uniformly,
  \item \label{item:14} $T$ is uniformly convex renormable,
  \item \label{item:19} $[L_2,T]$ is uniformly convex renormable,
  \item \label{item:16} $T$ is uniformly smooth renormable.
  \end{enumerate}
\end{theorem}
\begin{proof}
The equivalence of~\ref{item:15}--\ref{item:18} was shown in
Proposition~\ref{p:4}.

The equivalence of~\ref{item:11}--\ref{item:14} was shown in
Proposition~\ref{prop:5}.

Proposition~\ref{p:3} shows that~\ref{item:14} and~\ref{item:19} are
equivalent, hence all the properties~\ref{item:11}--\ref{item:14} for
$T$ and~\ref{item:11}--\ref{item:14} for $[L_2,T]$ are equivalent.

Assume that~\ref{item:15} does not hold. Then for some $c>0$ there are
infinitely many numbers $n$, such that
\[ \frac{\t(T|\os M_n)}{\sqrt n} \geq c.
\]
It follows from Lemma~\ref{lem:3} that for these numbers
\[ \frac{6\sqrt n}{\t^\circ(T|\os M_{2n})} \leq \frac{6\sqrt
  n}{\t(T|\os M_{2n})} \leq \frac 6c
\]
and therefore, by Proposition~\ref{prop:3} there exist factorizations
$\Sigma_n=B_n[L_2,T]A_n$ such that $\|B_n\|\,\|A_n\|\leq6/c$. Of
course, this yields such factorizations for \emph{all} $n\in\N$ and
therefore~\ref{item:12} cannot hold for $[L_2,T]$.

Hence we have shown that~\ref{item:12} implies~\ref{item:15}.

On the other hand Example~\ref{ex:3} shows that~\ref{item:18}
implies~\ref{item:12}.

So far, we have shown that~\ref{item:15}--\ref{item:19} are
equivalent.

Finally, since by Proposition~\ref{l:9}
properties~\ref{item:15}--\ref{item:18} for $T$ are equivalent to the
same properties for $T'$, it follows from
Proposition~\ref{p-convexity-dual} that also~\ref{item:16} is
equivalent to~\ref{item:14} and hence to all other properties.
\end{proof}

\begin{remark}
  Note that Pisier's Theorem~\ref{theorem:2} can easily be obtained
  from our Theorem~\ref{theorem:1} using the submultiplicativity of
  the ideal norms $\t(I_X|\os M_n)$. Namely, it follows from
  $\t(I_X|\os M_n)=o(n^{1/2})$ and the submultiplicativity, that there
  is a constant $c$ and a number $p_0>1$ such that
  \[ \t(I_X|\os M_n) \leq cn^{1/p_0-1/2} \quad\mbox{for all $n\in
    \N $.}
  \]
  By Pisier \cite{pis75}, this implies that $X$ has martingale type
  $p$ for all $p<p_0$.
\end{remark}

To formulate the last example, we will need one more definition.

\begin{definition*}
  For $T:X\to Y$, the \emph{Haar type $p$ ideal norm} $\t_p(T|\os
  H(\D_1^n))$, is the smallest constant $c\geq0$ such that
  \[ \Big\|\sum_{\D_1^n} [L_2,T]x_k^{(j)}\chi_k^{(j)} \Big|L_p\Big\|
  \leq c\, \Big(\sum_{k=1}^n \Big\|\sum_{j=1}^{2^{k-1}}
  x_k^{(j)}\chi_k^{(j)} \Big| L_p \Big\|^p\Big)^{1/p}
  \]
  for all $(x_k^{(j)})\subseteq X$.
\end{definition*}
In particular, $\t_2(T|\os H(\D_1^n))=\t(T|\os H(\D_1^n))$. An
operator $T$ is said to have martingale type $p$ if the sequence
$(\t_p(T|\os H(\D_1^n)))$ is bounded.

The following example shows that an analogue of the
Davis--Figiel--Johnson--Pe\l czy\'nski Theorem on weakly compact
operators (they factor through a reflexive Banach space; see
\cite[Cor.~1, p.~314]{davis74:_factor}) cannot hold for super weakly
compact operators, since if $T$ factors through a superreflexive
Banach space then it follows already that it has martingale type $p$
for some $p>1$; see \cite[Thm.~3.2, p.~340]{pis75}.

Let $t=(\tau_n)$ be a non-increasing sequence of positive numbers.  We
consider the diagonal operator $D_t:l_1\to l_1$ defined by
$D_t(\xi_k):=(\tau_k\xi_k)$.

The following fact is proved in \cite{wenzel:_haar}.
\begin{example}
  \label{ex:2}
  $\quad \t_p(D_t|\os H(\D_1^n)) = \Big( \sum_{k=1}^n
  |\tau_k|^{p'}\Big)^{1/p'}.$
\end{example}
\begin{corollary*}
  \label{cor:1}
  If $\tau_k=1/(1+\log k)$ then the operator $D_t$ is super weakly
  compact but does not have martingale type $p$ for any $p>1$ and
  hence does not factor through a superreflexive Banach space.
\end{corollary*}


\begin{thebibliography}{10}

\bibitem{beauzamy76:_op}
{\sc Beauzamy, B.}
\newblock Op\'erateurs uniform\'ement convexifiants.
\newblock {\em Studia Math.} {\bf 57} (1976), 103--139.

\bibitem{bea85a}
{\sc Beauzamy, B.}
\newblock {\em Introduction to {Banach} spaces and their geometry},
  volume~68 of {\em North-Holland mathematics studies}.
\newblock North-Holland, second edition, 1985.

\bibitem{bourgain89:_geomet_l}
{\sc Bourgain, J.}{\sc , Kalton, N.}{\sc ,\ {\rm and} Tzafriri, L.}
\newblock Geometry of finite dimensional subspaces and quotients of
  {$L_p$}.
\newblock In {\em Geometric aspects of functional analysis, Israel
  1987--88}, volume 1376 of {\em Lect. Notes in Math.}, pages
  138--175, 1989.

\bibitem{davis74:_factor}
{\sc Davis, W.}{\sc , Figiel, T.}{\sc , Johnson, W.~B.}{\sc ,\ {\rm
    and} {Pe\l czy\'nski}, A.}
\newblock Factoring weakly compact operators.
\newblock {\em J. Func. Anal.} {\bf 17} (1974), 311--327.

\bibitem{day41:_some}
{\sc Day, M.~M.}
\newblock Some more uniformly convex spaces.
\newblock {\em Bull. Amer. Math. Soc.} {\bf 47} (1941), 504--507.

\bibitem{enf72}
{\sc Enflo, P.}
\newblock Banach spaces which can be given an equivalent uniformly
  convex norm.
\newblock {\em Isr. J. Math.} {\bf 13} no.~3--4 (1972), 281--288.

\bibitem{heinrich80:_finit}
{\sc Heinrich, S.}
\newblock Finite representability and super-ideals of operators.
\newblock {\em Diss. Math.} {\bf 172} (1980), 37 p.

\bibitem{heinrich80:_ultrap_banac}
{\sc Heinrich, S.}
\newblock Ultraproducts in {Banach} space theory.
\newblock {\em J. Reine Angew. Math.} {\bf 313} (1980), 72--104.

\bibitem{hinrichs:_operat_radem_gauss}
{\sc Hinrichs, A.}
\newblock Operators of {Rademacher} and {Gaussian} subcotype.
\newblock Forschungsergebnisse der FSU Jena.

\bibitem{james72:_super_banac}
{\sc James, R.~C.}
\newblock Super-reflexive {Banach} spaces.
\newblock {\em Can. J. Math.} {\bf 24} no.~5 (1972), 896--904.

\bibitem{lindenstrauss63:_banac}
{\sc Lindenstrauss, J.}
\newblock On the modulus of smoothness and divergent series in
  {Banach} spaces.
\newblock {\em Mich. Math. J.} {\bf 10} (1963), 241--252.

\bibitem{lindenstrauss68:_absol_l}
{\sc Lindenstrauss, J.}{\sc \ {\rm and} {Pe\l czy\'nski}, A.}
\newblock Absolutely summing operators in {$\mathcal{L}_p$}-spaces and
  applications.
\newblock {\em Studia Math.} {\bf 29} (1968), 275--326.

\bibitem{pietsch:_orthon_banac}
{\sc Pietsch, A.}{\sc \ {\rm and} Wenzel, J.}
\newblock Orthonormal systems and {Banach} space geometry.
\newblock Cambridge University Press, 1998.

\bibitem{pis75}
{\sc Pisier, G.}
\newblock Martingales with values in uniformly convex spaces.
\newblock {\em Isr. J. Math.} {\bf 20} no.~3--4 (1975), 326--350.

\bibitem{wenzel:_haar}
{\sc Wenzel, J.}
\newblock Haar type ideal norms of diagonal operators.
\newblock preprint, available from
  {\url|http://www.minet.uni-jena.de/~wenzel/example.tex|}.

\end{thebibliography}

\end{document}